\documentstyle{amsppt}
\magnification=1200
\catcode`\@=11
\redefine\logo@{}
\catcode`\@=13

\define \bn{\Bbb N}
\define \bz{\Bbb Z}
\define \bq{\Bbb Q}
\define \br{\Bbb R}
\define \bc{\Bbb C}

\define \M{{\Cal M}}
\define\Ha{{\Cal H}}
\define\La{{\Cal L}}

\define\rk{\text{rk}~}


\define\NEF{\text{NEF}}

\define\NE{\text{NE}}
\define\NS{\text{NS}}
\define\Exc{\text{Exc}}
\TagsOnRight

\document

\topmatter
\title
A remark on algebraic surfaces with polyhedral Mori cone
\endtitle

\author
Viacheslav V. Nikulin \footnote{Supported by
Grant of Russian Fund of Fundamental Research
\hfill\hfill}
\endauthor

\address
Steklov Mathematical Institute,
ul. Gubkina 8, Moscow 117966, GSP-1, Russia
\endaddress
\email
slava\@nikulin.mian.su
\endemail

\dedicatory
To 75-th Birthday of I.R. Shafarevich
\enddedicatory

\abstract
We denote by $FPMC$ the class of all non-singular projective
algebraic surfaces $X$ over $\bc$ with finite polyhedral Mori cone
$\NE(X)$ $\subset \NS(X)\otimes \br$. If $\rho(X)=\rk \NS(X)\ge 3$, then
the set $\Exc(X)$ of all exceptional curves on $X\in FPMC$ is finite and
generates $\NE(X)$. Let $\delta_E(X)$ be the maximum of $(-C^2)$ and
$p_E(X)$ the maximum of $p_a(C)$ respectively for all $C\in \Exc(X)$.
For fixed $\rho\ge 3$, $\delta_E$ and $p_E$ we denote by
$FPMC_{\rho,\delta_E,p_E}$ the class of all
algebraic surfaces $X\in FPMC$ such that
$\rho(X)=\rho$, $\delta_E(X)=\delta_E$
and $p_E(X)=p_E$. We prove that the class
$FPMC_{\rho,\delta_E,p_E}$
is bounded in the following sense:
for any $X\in FPMC_{\rho,\delta_E,p_E}$
there exist an ample effective divisor $h$ and a very ample divisor
$h'$ such that $h^2\le N(\rho,\,\delta_E)$ and
${h'}^2\le N'(\rho,\,\delta_E,\,p_E)$
where the constants $N(\rho,\,\delta_E)$ and $N'(\rho,\,\delta_E,\,p_E)$
depend only on $\rho,\,\delta_E$ and $\rho,\,\delta_E,\,p_E$
respectively.

One can consider Theory of surfaces $X\in FPMC$ as Algebraic Geometry
analog of the Theory of arithmetic reflection groups
in hyperbolic spaces.
\endabstract

\rightheadtext
{Algebraic surfaces and Mori cone}
\leftheadtext{V.V. Nikulin}

\endtopmatter

\head
0. Introduction
\endhead

Let $X$ be a non-singular projective algebraic surface over
algebraically closed field with finite polyhedral Mori cone
$\NE(X)\subset \NS(X)\otimes \br$ where $\NS(X)$ is the Neron--Severi lattice
of $X$. If $\rho=\rk \NS(X)\ge 3$, then the set
$\Exc(X)$ of exceptional curves of $X$ is finite and generates the
cone $\NE(X)$. Further we assume that $\rho \ge 3$.
One can introduce natural invariants of $X$:
$$
\rho=\rk \NS(X),\ \ \delta_E=\max_{C\in \Exc(X)}{(-C^2)},\ \
p_E=\max_{C\in \Exc(X)}{p_a(C)}.
$$
The main result of the paper (Theorem 1.1)
is that {\it the class of surfaces $X$ with finite polyhedral
Mori cone and fixed invariants $\rho\ge 3$, $\delta_E$ and
$p_E$ is bounded:} there exists an effective ample
divisor $h$ and a very ample divisor $h^\prime$ on $X$ such that
$h^2\le N(\rho,\,\delta_E)$ and ${h^\prime}^2\le
N(\rho,\,\delta_E,\,p_E)$.

The key step in the proof the theorem is using
the old result of the author \cite{N4}
on ``narrow parts'' of convex finite polyhedra of finite volume in
hyperbolic spaces (for the special case of Theorem 1.1 it is
formulated in Lemma 1.1). It was used in \cite{N4}
to prove some finiteness results on arithmetic reflection
groups in hyperbolic spaces.
We also use some standard results about
symmetric matrices with non-negative coefficients
(Perron--Frobenius Theorem) and Reider's Theorem \cite{R}
on very ample divisors of surfaces. These considerations permit to
write down the $h$ and $h^\prime$ as
linear combination of exceptional curves on the $X$.

In Example 1.2 we show that Theorem 1.1 is not valid if
one of the invariants $\rho\ge 3$, $\delta_E\ge 3$ and $p_E$ is not fixed.

Because of Theorem 1.1, one can ask about classification of surfaces
$X$ with finite polyhedral Mori cone and small invariants
$\rho$, $\delta_E$ and $p_E$.

In Example 1.3 we consider classification for $\delta_E=1$. Using
results of \cite{N2} and \cite{N3},
we then have $\rho\le 9$. In Example 1.3.1 we additionally to
$\delta_E=1$ suppose that $p_E=0$. Then
one gets non-singular Del Pezzo surfaces whose classification
is well-known, e. g. see \cite{Ma}.

In Example 1.4 we consider classification for $\delta_E=2$.
Using results of \cite{V1} and \cite{E}, we have $\rho\le 22$.
In Example 1.4.1 we additionally to $\delta_E=2$
suppose that $p_E=0$. Then one gets (minimal) K3 surfaces,
(minimal) Enriques surfaces, minimal resolutions of
singularities of Del Pezzo surfaces with Du Val
singularities, and rational surfaces with $K^2=0$ and nef $-K$.
Classification of the last class of surfaces with finite
polyhedral Mori cone was not considered in literature, and we
give this classification. It uses Ogg--Shafarevich theory of
elliptic surfaces and results \cite{Ha}, \cite{D} and \cite{CD}
about rational elliptic surfaces. It is interesting that not all
rational surfaces with $K^2=0$, nef $-K$ and finite polyhedral
Mori cone are elliptic.

It seems, nobody tried to classify surfaces $X$ with finite polyhedral
Mori cone and $p_E\ge 1$.

It Sect. 2 we consider generalization of the main theorem 1.1 above to
some surfaces with locally finite polyhedral Mori cone.

I am grateful to I.V. Dolgachev and I.R. Shafarevich for very useful
discussion on elliptic surfaces.

\head
1. Algebraic surfaces with finite polyhedral Mori cone
\endhead

Let $X$ be a non-singular projective
algebraic surface over an algebraically closed field. Let $\NS(X)$
be {\it Neron--Severi lattice} of $X$
(i. e. the group of divisors on $X$ by numerical equivalence
considered together with the intersection pairing). By Hodge Index Theorem,
the lattice $\NS(X)$ is hyperbolic: it has signature $(1, \rho-1)$ where
$\rho=\rk \NS(X)$. We denote by
$\NE(X)\subset \NS(X)\otimes \br$ the {\it Mori cone} of $X$
generated over $\br_+=\{t\in \br\,|\,t\ge 0\}$
by all effective curves on $X$.
By definition, a surface $X$ has a {\it finite polyhedral Mori cone}
$\NE(X)$ if $\NE(X)$ is generated by a finite set of rays (we denote by
$FPMC$ the class of all these surfaces).
The minimal set of these rays is called the
{\it set of extremal rays}. We denote by
$V^+(X)$ the half-cone containing a polarization (i. e.
an ample divisor) of the light cone
$V(X)=\{x\in \NS(X)\otimes \br\ \mid \ x^2>0\}$. By Riemann--Roch
Theorem, the cone $\NE(X)$ contains the half-cone $V^+(X)$.
It follows that for $\rho(X)\ge 3$,
the set of extremal rays of $X\in FPMC$ is equal to
$\br_+E$, $E\in \Exc(X)$, where $\Exc(X)$ is the set of all
exceptional (i. e. irreducible and having negative square) curves of $X$.
In particular, the set $\Exc(X)$ is finite. Thus, we can introduce
{\it natural invariants} of $X\in FPMC$:
$$
\rho(X)=\rk \NS(X),
\tag{1.1}
$$
$$
\delta_E(X)=\max_{C \in \Exc(X)}{(-C^2)},
\tag{1.2}
$$
$$
p_E(X)=\max_{C \in \Exc(X)}{p_a(C)},
\tag{1.3}
$$
where $p_a(C)={C^2+C\cdot K\over 2}+1$ is the arithmetic genus of
a curve $C$, $K$ is the canonical divisor of $X$.

Surfaces $X\in FPMC$ are interesting because of the following reasons:

1) Polyhedrality of the Mori cone $\NE(X)$ is very important in
Mori Theory (see \cite{Mo}).
It is interesting and curious to ask what will
be if one requires the only this condition.

2) We consider surfaces $X\in FPMC$ as
Algebraic Geometry analog  of arithmetic groups generated by reflections in
hyperbolic spaces (e. g. see \cite{N4}, \cite{N5} and \cite{N8}).
We also expect that they are connected with some analog  of
automorphic products introduced by R. Borcherds
(see \cite{B1}, \cite{B2}, \cite{GN1}---\cite{GN7} and
\cite{N10}---\cite{N12}).

3) We expect that quantum cohomology related with
surfaces $X\in FPMC$ are very interesting: one can consider the
set $\Exc (X)$ as analog of a system of simple real roots.
Here ``related'' means that not necessarily the quantum cohomology of
$X$ itself, but e. g. quantum cohomology of varieties fibrated by
$X\in FPMC$ might be interesting ones.
See some examples in \cite{CCL}, \cite{HM},
\cite{Ka}, \cite{Moo}  and also \cite{GN3}, \cite{GN7}.

\remark{Basic example 1.1} We have the following basic
example of surfaces $X\in FPMC$.
It shows that there are plenty of
surfaces $X\in FPMC$. For this example, let $Y$ be
a normal projective algebraic
surface such that $-K_Y$ is nef and $K_Y^2>0$ (in particular, one
can take $Y$ to be a numerical
Del Pezzo surface with normal singularities).
Let $X$ be the minimal
resolution of singularities of $Y$. Then $X\in FPMC$ if
$\rho (X)\ge 3$.
This follows from the Mori Theory \cite{Mo} applied to
the non-singular projective algebraic surface $X$.
E. g. see \cite{N9}.
\endremark

\smallpagebreak

For fixed invariants $\rho\ge 3$, $\delta_E$, $p_E$,
we denote by $FPMC_{\rho,\delta_E,p_E}$ the class of all
algebraic surfaces $X\in FPMC$ such that
$\rho(X)=\rho$, $\delta_E(X)=\delta_E$ and $p_E(X)=p_E$.

In Theorem 1.1 below we want to show that the class
$FPMC_{\rho,\delta_E,p_E}$ is bounded.
We remind that any non-singular projective
algebraic surface has a linear projection embedding into
${\Bbb P}^5$ (e. g. see \cite{Sh2}). This projection
keeps the degree. Surfaces in ${\Bbb P}^5$ of the fixed degree
depend on a finite number of Chow coordinates
(e. g. see \cite{Sh2}).

\proclaim{Theorem 1.1} For $\rho\ge 3$, there are
constants $N(\rho,\,\delta_E)$ and $N'(\rho,\,\delta_E,\,p_E)$ depending
only on $(\rho,\,\delta_E)$ and $(\rho,\,\delta_E,\,p_E)$ respectively
such that for any $X\in $ \linebreak
$FPMC_{\rho,\delta_E,p_E}$ there exists an
ample effective divisor $h$ such that $h^2\le N(\rho,\,\delta_E)$, and
if the ground field is $\bc$,
there exists a very ample divisor $h'$ such that
${h'}^2\le N'(\rho,\,\delta_E,\,p_E)$.
\endproclaim

\demo{Proof} Let $\NEF(X)=\NE(X)^{\,\ast}$ be the dual nef cone.
We have
$\NEF(X)\subset\overline{V^+(X)}$ $=\overline{V^+(X)}^{\ \ast}\subset \NE(X)$.
Therefore, the nef cone defines
a finite polyhedron $\M=\NEF(X)/\br_+$ of finite volume
in the hyperbolic space $\La(X)=V^+(X)/\br_{++}$, where
$\br_{++}=\{t \in \br\ |\ t>0\}$. The set
$\Exc(X)$ is the set of orthogonal vectors to faces (of the highest
dimension) of $\M$. By \cite{N4, Appendix, Theorem 1} (see
also \cite{N15} about much more exact statements), we have
the following

\proclaim{Lemma 1.1} There exist exceptional curves
$E_1,\dots ,E_\rho\in \Exc(X)$ such that the conditions
(a), (b) and (c) below are valid:

(a) $E_1,\dots ,E_\rho$ generate $\NS(X)\otimes\bq$;

(b) ${2\,(E_i\cdot E_j)\over \sqrt{E_i^2E_j^2}}<62$;

(c) The dual graph of $E_1,\dots ,E_\rho$ is connected, i. e. one cannot
divide this set in two non-empty subsets orthogonal to one another.
\endproclaim

Let us consider the $\rho \times \rho$
matrix $\Gamma=(\gamma_{ij})=(E_i\cdot E_j)$ where
$1\le i,\,j\le \rho$. By conditions of the Theorem 1.1 and
by Lemma 1.1, we have
$-\delta_E \le \gamma_{ii}<0$ for any $1\le i\le \rho$,
and $0\le \gamma_{ij} < 31\delta_E$ for any $1\le i,\,j\le \rho$ and
$i\not=j$. It follows that the set of possible matrices $\Gamma$ is
finite. Thus, in further considerations
we can fix one of the possible matrices $\Gamma$.

The matrix $\Gamma$ has non-negative
coefficients out of its diagonal. Moreover, it is symmetric and
indecomposable (by the condition (c)).
Thus, by Perron-Frobenius Theorem,
its maximal eigenvalue
$\lambda$ has
multiplicity one and has the eigenvector
$v=b_1E_1+ \cdots + b_\rho E_\rho$ with positive coordinates $b_i>0$.
Since the lattice $\NS(X)$ is hyperbolic, the eigenvalue $\lambda >0$.
It follows that
$\Gamma v=\lambda v$ and $E_j\cdot v=\lambda b_j>0$ for any
$1\le j\le \rho$. Moreover, $v^2=\lambda(b_1^2+\cdots +b_\rho^2)>0$.

We can replace real numbers $b_i$ by very closed positive
rational numbers ${b'}_i$ keeping the inequalities
$E_j\cdot ({b'}_1E_1+\cdots +{b'}_\rho E_\rho)>0$
for any $1\le j\le \rho$, and
$({b'}_1E_1+\cdots +{b'}_\rho E_\rho)^2>0$.
Multiplying numbers ${b'}_i$ by an appropriate positive natural
number $N$, finally we find natural numbers $a_i=N{b'}_i$ such that
$E_j\cdot (a_1E_1+\cdots +a_\rho E_\rho)>0$
for any $1\le j\le \rho$, and $(a_1E_1+ \cdots +a_\rho E_\rho)^2>0$.
It follows

\proclaim{Lemma 1.2} Under the conditions (a) and (c) of
Lemma 1.1, there exist $a_i\in \bn$, $i=1, \dots ,\rho$,
depending only on the matrix $\Gamma=(E_i\cdot E_j)$,
$1\le i,\,j\le \rho$, such that for
$h=a_1E_1+ \cdots +a_\rho E_\rho$
one has
$$
E_j\cdot h>0
$$
for any $1\le j\le \rho$, and $h^2>0$.
\endproclaim

Suppose that $C$ is an irreducible curve on $X$ different from
$E_1,\dots, E_\rho$. Then $C$ defines a non-zero element in
$\NS(X)$ because $C\cdot H>0$ for a hyperplane section $H$.
It follows that $C\cdot E_i\ge 0$
for any $0\le i\le \rho$, and at least one of these inequalities is strong
because $E_1,\dots, E_\rho$ generate $\NS(X)\otimes \bq$. Thus,
$C\cdot h>0$. It follows that $C\cdot h>0$ for any effective curve $C$.
Since $h^2>0$,
by Nakai---Moishezon criterion, the divisor $h$ is ample. By
the construction, $h$ is effective.
Since for the fixed $(\rho,\,\delta_E)$
the set of possible matrices $\Gamma$ is finite,
$h^2\le N(\rho,\,\delta_E)$ for a constant $N(\rho,\,\delta_E)$ depending only
on $(\rho,\,\delta_E)$. It proves the first statement of the Theorem 1.1.

To prove second statement, let us additionally
fix $p_a(E_i)$, $1\le i\le \rho$. Since $0\le p_a(E_i)\le p_E$,
there exists only a finite number of possibilities.
Let $K=d_1E_1+\cdots+d_\rho E_\rho$ be the canonical class of $X$
where $d_i\in \bq$.
One can find $d_1,\dots , d_\rho$ from
equations ${E_i^2+(E_i\cdot K)\over 2}+1=p_a(E_i)$, $1\le i\le \rho$.
Since the lattice $\NS(X)$ is non-degenerate and $E_1,\dots ,E_\rho$
give a bases of $\NS(X)\otimes \bq$, these equations define
$d_1,\dots,d_\rho$ uniquely.

Suppose that the ground field is $\bc$.
By Reider's Theorem \cite{R} (see also \cite{L}),
we have $h^\prime=K+4h$ is very ample.
It follows that ${h^\prime}^2\le N'(\rho,\,\delta_E,\,p_E)$
where $N'(\rho,\,\delta_E,\,p_E)$ is bounded by a constant depending on
$(\rho,\,\delta_E,\,p_E)$. (We remark that this is the only place in
the proof of Theorem 1.1 where we use that the ground
field is $\bc$.) This finishes the proof of Theorem 1.1.
\enddemo

In Theorem 1.1, we fixed invariants $\rho$, $\delta_E$ and
$p_E$. In example 1.2 below, we show that the theorem 1.1 is not
in general true if one does not fix one of these invariants.

\remark{Example 1.2} Let us consider a non-singular curve $C_g$ of
genus $g$ and an invertible sheaf $\Cal L$ on $C_g$
of the degree $-n$ where $n\in \bn$. The ruled surface
$\pi: Y={\Bbb P}({\Cal O}_{C_g}\oplus {\Cal L})\to C_g$
has the exceptional section $C_g$ of genus $g$ such that
$(C_g)^2=-n$ (e. g. see \cite{Har, Ch. V, Example 2.11.3}).
Let $E_0$ be a fiber of $\pi$ and
$X$ the blow up of $Y$ in a point of $E_0$ which does not
belong to the section $C_g$. Let $F_0$ be the exceptional curve
of the blow up. We get three exceptional curves
$C_g$, $E_0$ and $F_0$ on $X$ of the genus $g$, $0$ and $0$
respectively and with the intersection matrix
$$
\pmatrix
-n&1&0\\
1 &-1&1\\
0 & 1&-1
\endpmatrix.
$$
Using this matrix, it is easy to prove that the cone
$\br_+C_g+\br_+E_0+\br_+F_0$ contains the cone $V^+(X)$. (One needs to
show that the Gram matrix of any two of these three curves is
either negative or semi-negative definite.) It
follows that $\NE(X)=\br_+C_g+\br_+E_0+\br_+F_0$, and
$X\in FPMC_{3,n,g}$ where the numbers
$n>0$ and $g\ge 0$ can be arbitrary.
The surfaces $X$ give the infinite dimensional family of surfaces
in $FPMC_{3,\delta_E=n}=\cup_{p_E\ge 0}{FPMC_{3,n,p_E}}$. This shows that
Theorem 1.1 is not true if one does not fix the invariant $p_E$.
The same example shows that Theorem 1.1 is not true if one does
not fix the invariant $\delta_E$.

Now let us show that $\rho$ is not bounded for
$X\in FPMC_{\delta_E,p_E}=$ \linebreak
$\cup_{\rho\ge 3}{FPMC_{\rho,\delta_E,p_E}}$
for any fixed $\delta_E\ge 3$ and $p_E\ge 0$.
Let us consider the surface $X$ above with
$n=\delta_E\ge 3$ and $g=p_E$.
It has three exceptional curves
$C_g$, $E_0$ and $F_0$. Curves $E_0$ and $F_0$ (different from
$C_g$) define a connected tree of curves and have
one intersection point of the curves.
Consider the blow up $X_1$ of $X$ in the
intersection point. Then $X_1$ has four exceptional curves
$C_g$, $E_0$, $F_0$ and $F_1$ where $F_1$ is the exceptional
curve of the blow-up. Like for Example 1.1, one can show that
$X_1\in FPMC$ and $C_g$, $E_0$, $F_0$ and $F_1$ are all exceptional
curves of $X_1$. It has $\rho(X_1)=4$. Curves $E_0$, $F_0$ and
$F_1$ (different from $C_g$) also define a tree
and have two intersection points. We can repeat this procedure
considering blow up $X_2$ of $X_1$ in one of these two points.
Repeating this procedure, we get an infinite sequence $X_k$, $k\ge 0$,
of surfaces with $\rho(X_k)=3+k$. We have: $X_k\in FPMC$ and $X_k$
has exactly $3+k$ exceptional curves where $2+k$ of them are proper
preimages of exceptional curves of $X_{k-1}$ and one is the
exceptional curve of the blow-up in an intersection point of
two exceptional curves (different from $C_g$) of $X_{k-1}$.
One can see that this sequence contains an infinite sequence of surfaces
$X_k\in FPMC_{\rho=3+k,\delta_E,p_E}$.
Actually one can find an infinite
sequence $X_k$ such that the curve $C_g$ of $X_k$
has $(C_g)^2=-n=-\delta_E$, and all other exceptional curves
$E$ of $X_k$ have $-3\le E^2<0$. It follows, that
$X_k\in FPMC_{\delta_E,p_E}$ since we assume that $\delta_E\ge 3$.
\endremark

\smallpagebreak

Because of Theorem 1.1, one can ask about classification of
surfaces \linebreak
$X\in FPMC_{\rho,\delta_E,p_E}$ for small invariants
$\rho$, $\delta_E$ and $p_E$.

\smallpagebreak

\remark{Example 1.3} For this example, we suppose that
$\delta_E=1$. Then $\rho\le 9$. Really, for $\delta_E=1$
any exceptional curve $E\in \Exc(X)$ has $E^2=-1$ and defines a reflection
of $\NS(X)$ which maps $E\to -E$ and is identical on $E^\perp$.
It is given by the formula $x\mapsto x+2(E\cdot x)E$, $x\in \NS(X)$.
All $E\in \Exc(X)$ generate a reflection group $W\subset O(\NS(X))$ with
the fundamental chamber $\NEF(X)/\br_{+}\subset \La (X)$
of finite volume. It follows that
$[O(\NS(X)):W^{(-1)}(\NS(X))]<\infty$ where
$W^{(-1)}(\NS(X))$ is generated by reflections in all $\alpha\in \NS(X)$
with $\alpha^2=-1$.
It was shown in \cite{N2}, \cite{N3} that
$\rk S \le 9$ for any hyperbolic lattices $S$ with
$[O(S):W^{-1}(S)]<\infty$, and all hyperbolic lattices $S$ with
this property were found; see also \cite{N6}, \cite{N8}.
From Theorem 1.1, it then follows that
the family $FPMC_{\delta_E=1,p_E}$ is bounded and may be
described (in principle) for a fixed $p_E$. Here we denote
$FPMC_{\delta_E=k,p_E}=\cup_{\rho\ge 3}{FPMC_{\rho,k,p_E}}$.
\endremark

\smallpagebreak

\remark{Example 1.3.1} Let us additionally (to the
condition $\delta_E=1$) assume that $p_E=0$.
We have: {\it The family $X\in FPMC_{\delta_E=1,p_E=0}$
consists of all non-singular
Del Pezzo surfaces $X$ and is well-known} (e. g. see \cite{Ma}).
Really,
we have that $-K\cdot E=1>0$ for any $E\in \Exc(X)$. Since $\Exc(X)$
is finite and generates $\NE(X)$, it then follows that $(-K)^2>0$, and
$-K$ is ample by Nakai---Moishezon criterium. The opposite statements
is a very particular case of Basic Example 1.1. One can get Del Pezzo
surfaces as blow up of ${\Bbb P}^2$ in $\le 8$ points in ``general''
position. It follows that they define a bounded family of algebraic
surfaces and illustrates the Theorem 1.1 for this very particular
case. We even have more: their moduli have finite number of connected
components.

It seems, nobody tried to classify $FPMC_{\delta_E=1,\,p_E}$
for $p_E\ge 1$.
\endremark

\smallpagebreak

\remark{Example 1.4} For this example, we suppose that
$\delta_E=2$. Thus, $E^2=-1$ or $-2$ for any $E\in \Exc(X)$.
Then again any $E\in \Exc(X)$ defines
reflection of $\NS(X)$. It is given by the formula
$x\mapsto x-\left(2(E\cdot x)/E^2\right)E$, $x\in \NS(X)$.
The same arguments as above show that
$\rho\le 22$; see \cite{V1}, \cite{V2} and \cite{E}. We mention that
the same result is valid for $X\in FPMC$
if all $E\in \Exc(X)$ define reflections of $\NS(X)$.
We also remark that
if $E^2=-2$ for any $E\in \Exc(X)$, then $\rho\le 19$,
see \cite{N2}, \cite{N3}.
By Theorem 1.1, it then follows that
the class $FPMC_{\delta_E=2,\,p_E}$ is bounded for a fixed $p_E$.
\endremark

\smallpagebreak

\remark{Example 1.4.1} Let us additionally (to the
condition $\delta_E=2$) assume that $p_E=0$.
Let $X\in FPMC_{\delta_E=2,\,p_E=0}$.
Then $E^2=-1$ or $-2$ and $E$ is non-singular rational for any
$E\in \Exc(X)$.
By the formula for genus of curve, $-K\cdot E=1$ if $E^2=-1$,
and $-K\cdot E=0$ if $E^2=-2$. It follows that $-K$ is nef and
$K^2\ge 0$. Vice versa, by the formula for genus of curve,
any surface $X$ with nef $-K$ has only non-singular rational
exceptional curves $E$ with $E^2=-1$ or $-2$.
Considering these surfaces $X\in FPMC$, we get
one of cases below where we for simplicity suppose that the
basic field $k=\bc$ (one can consider arbitrary $k$ using results
from \cite{CD}).

{\bf Case 1.} {\it Suppose that $E^2=-2$ for any
$E\in \Exc(X)$.} Then $K\equiv 0$ and $X$ is minimal.
By classification of surfaces
(e. g. see \cite{Sh1}),
we then have $mK=0$ for some $m\in \bn$.
We can suppose that $m$ is minimal with
this property. Then $K$ defines the cyclic $m$-sheeted covering
$\pi:\widetilde{X}\to X$ where $K_{\widetilde{X}}=0$ and
$\widetilde{X}$ is either $K3$ or Abelian surface. The preimage
$\pi^{-1}(E)$ of an exceptional curve $E$ of $X$ contains an exceptional
curve of $\widetilde{X}$ because $\pi^{-1}(E)^2=mE^2<0$.
Abelian surfaces do not have exceptional curves.
Thus, $\widetilde{X}$ is K3 surface and $m=1$ or $2$.
It follows that $X$ is either K3 (if $m=1$), or Enriques
(if $m=2$) surface. Thus {\it we get K3 or Enriques
surfaces $X$.} All K3 and Enriques surfaces with
finite polyhedral Mori cone (equivalently, when their automorphism
group is finite \cite{ P-\u S\u S})
were classified in \cite{N2}, \cite{N3}, \cite{N6}, \cite{N7},
\cite{N8} and \cite{Ko}. Their moduli have finite number of connected
components. It is interesting that to get this classification,
general arguments of Lemma 1.1 were used in some cases (see \cite{N6}).

{\bf Case 2.} {\it Suppose that there exists $E\in \Exc(X)$ with
$E^2=-1$.} Then $K\cdot E=-1$ and $K$ is not numerically zero.
Since $-K$ is nef, it follows that linear systems
$|nK|$ are empty for all $n>0$.
By classification of surfaces, $X$ is either ruled or rational.
Suppose that $X$ is ruled and $\pi:X\to C$ is a morphism on a curve $C$.
Since $\Exc(X)$ generates $\NS(X)\otimes \bq$,
the image of one of exceptional curves $E$ is equal to $C$.
Since $E$ is rational, the curve $C$ is rational. It follows that
{\it $X$ is a rational surface with nef $-K$.}

{\bf Case 2a.} {\it Suppose that $X$ is rational,
$-K$ is nef and $K^2>0$.}
Since $E\cdot K=0$ for any $E\in \Exc(X)$ with $E^2=-2$, the intersection
matrix of these curves is negative definite. It follows that
these curves define several disjoint configurations of the
Dynkin type $A_n$, $D_n$ or $E_n$. It follows that there exists
the contraction map $\pi: X \to Y$ of these curves where
$Y$ is a surface with Du Val singularities. We then have
$K_X=\pi^\ast K_Y$. It follows that $(-K_Y)$ is numerically ample,
and $Y$ is Del Pezzo surface with Du Val singularities. It is a
particular case of Basic Example 1.1.
Thus, we get that {\it the surfaces $X$ are minimal resolutions of
singularities of Del Pezzo surfaces with Du Val singularities.}
Classification of these surfaces
as blow-ups of ${\Bbb P}^2$ in $\le 8$ points in ``general'' position
was obtained by Nagata \cite{Na}. See also \cite{Ma}.
It follows that these surfaces define a bounded family of algebraic
surfaces and illustrates the Theorem 1.1 for this very particular
case. Their moduli have finite number of connected
components. Classification of all possible graphs
$\Gamma (\Exc(X))$ of exceptional curves on $X$ was obtained in
\cite{AN1} and \cite{AN2}.

{\bf Case 2b.} {\it Suppose that $X$ is rational, $-K$ is nef
and $K^2=0$.}
By Noether formula, $\rho=\rk \NS(X)=10$.
Since $\NE(X)$ is finite polyhedral,
the dual cone $\NEF(X)\subset \overline{V^+(X)}$
is also finite polyhedral and the polyhedron $\M(X)=$
\linebreak $\NEF(X)/\br_{+}$ is a finite polyhedron in
the hyperbolic space $\La(X)$. The polyhedron $\M(X)$ has
$\br_{++}(-K)$ as its infinite vertex and should be finite
in a neighbourhood of this vertex. It follows that
the set of exceptional curves $E\in \Exc(X)$ which are orthogonal
to $-K$ is parabolic and has the rank 8. It means that any connected
component $\Gamma_i$ of the dual graph $\Gamma$ of the set of
these curves has semi-negative definite Gram matrix of the rank $r_i$ and
the sum $\sum_i{r_i}=8$. Since $E^2=-2$ if $E\in \Exc(X)$ and
$E\cdot K=0$, it follows that the graphs $\Gamma_i$ are extended
Dynkin diagrams of types $\widetilde{A}_{r_i}$, $\widetilde{D}_{r_i}$
or $\widetilde{E}_{r_i}$ where $\sum_i{r_i}=8$.

Let us show that the opposite statement is also true. Suppose that
$X$ is rational, $-K$ is nef, $K^2=0$, any connected component
of the set of exceptional curves $E\in \Exc(X)$ with $E^2=-2$ is
extended Dynkin diagram of the rank $r_i$ and $\sum_i{r_i}=8$.
We claim that then $X\in FPMC$.

Since $-K$ is nef, $K^2=0$ and $K$ is not numerically zero,
the polyhedron $\M(X)=\NEF(X)/\br_{+}\subset \La(X)$ is
parabolic relative to $\br_{++}(-K)$.
E. g. see \cite{N9} (this follows from the Mori Theory applied
to non-singular surfaces). This means that $\M(X)$ is
finite in any angle of the hyperbolic space $\La (X)$
with the vertex  $\br_{++}(-K)$. (This is example of surfaces with
almost finite polyhedral Mori cone which we shall consider in Sect. 2).
Let $Q$ be the set of all exceptional curves $E\in \Exc(X)$ with $E^2=-2$
and $\Ha_E^+=\{0\not= x\in \NS(X)\otimes \br\ |\ x\cdot E\ge 0\}/\br_{++}$.
We have $\M(X)\subset \bigcap_{E\in Q}\Ha_E^+$.
Curves $E\in Q$ are all orthogonal to $-K$, and by the
condition on the set of these curves, the set
$\bigcap_{E\in Q}\Ha_E^+$ is a finite polyhedral angle in $\La(X)$
 with the vertex $\br_{++}(-K)$. It follows that
$\M(X)$ is finite polyhedral in a neighbourhood of
$\br_{++}(-K)$ in $\La(X)$.
It follows that $\M(X)$ and $\NEF(X)$ are finite.
Thus, the dual cone $\NE(X)$ is also finite. This proves the statement.

\noindent
{\it (*) Thus, the surfaces $X$ are rational surfaces with
nef $-K$, $K^2=0$, and such that any connected component
of the set of exceptional curves $E\in \Exc(X)$ with $E^2=-2$ is
an extended Dynkin diagram $\Gamma_i$
of the rank $r_i$ and $\sum_i{r_i}=8$.}

Classification of surfaces from (*) was not considered in
literature. So, we are forced to give this classification
below.

Suppose that $X$ satisfies (*). Let $E_j$, $j=1,\dots, r_i+1$
are exceptional curves with $E_j^2=-2$
which give an extended Dynkin diagram $\Gamma_i$.
There are natural $a_j$, $j=1,\,\dots,\,r_i+1$
such that $D_i =\sum_j{a_jE_j}$ has
$D_i^2=0$ and one of coefficients $a_j$ is equal to one.
Obviously, $D_i\in |-m_iK|$ where $m_i$ is the invariant of
the divisor $D_i$ and of the connected component $\Gamma_i$ of
the set of exceptional curves with square $-2$ of $X$.
By Riemann--Roch Theorem,
the linear system $|-K|$ is not empty and contains a divisor $D$.
If $D$ is different from one of divisors $D_i$, the linear
system $|-m_iK|$ has positive dimension. It follows that for
some natural $m$ the linear system $|-mK|$ is a pencil. By
Bertini Theorem, this pencil contains a non-singular curve
which is an elliptic curve. Thus, $|-mK|$ is elliptic pencil, and
the surface $X$ is then rational minimal elliptic surface
(see \cite{D} and \cite{CD} about theory of rational elliptic
surfaces). We shall consider classification of elliptic surfaces from
(*) in more details below. Our consideration also shows that $X$
might be not elliptic only if there exists
exactly one extended Dynkin diagrams
$\Gamma_1$, it has rank 8 and the divisor $D_1\in |-K|$.

Thus, the set of surfaces from (*) and its classification is
divided in two parts: Elliptic surfaces from (*), and

\noindent
{\it (***) rational surfaces $X$ with nef $-K$ and $K^2=0$ such
that all exceptional curves $F_j$ with $F_j^2=-2$ of $X$ define a
connected extended Dynkin diagram of the rank $8$.
There exist natural $a_j$ such that one of $a_j$ is equal to one and
$D=\sum_j{a_jF_j}\in |-K|$.}

First we consider classification of surfaces (***). The
graph $\Gamma$ of $F_1,\dots , F_9$ has type $\widetilde{E_8}$,
$\widetilde{D_8}$ or $\widetilde{A_8}$. Since $X$ is not relatively
minimal, there exists an exceptional curve $E$ of the first kind.
We have $E\cdot D=1$. It follows that $E\cdot F_{j_0}=1$ for
one of curves $F_{j_0}$ having $a_{j_0}=1$, and $E\cdot F_j=0$ if
$j\not=j_0$.

Suppose that the graph $\Gamma=\widetilde{E_8}$. Then
the exceptional curves $F_j$ and $E$ have the dual graph
$H\widetilde{E_8}$ below where black vertices correspond to
the curves $F_j$ with $F_j^2=-2$. An edge of the graph means transversal
intersection of curves in one point.
Analyzing the graph $H\widetilde{E_8}$, one can see that
the dual cone
$\{x\in \NS(X)\otimes \br \ |\ x\cdot E_i\ge 0\}$
of curves $E_1,\dots, E_{10}$ is contained in $\overline{V^+(X)}$.
It follows that $X \in FPMC$, and the curves $E_1,\dots , E_{10}$
are all exceptional curves of $X$. Thus, $H\widetilde{E_8}$
is the graph of all exceptional curves on $X$ if
$\Gamma=\widetilde{E_8}$.

If $\Gamma=\widetilde{D_8}$, considerations are a little bit more
complicated. Curves $F_j$ and $E$ define a subgraph of the
graph $H\widetilde{D_8}$ below with vertices $E_1,\dots ,E_{10}$.
Curves $E_2$, $E_4$, $E_6$, $E_8$ and $E_9$ of the subgraph
define an extended Dynkin diagram of the type $\widetilde{B_4}$.
The divisor $C=E_2+E_4+2E_6+2E_8+2E_9$ is nef
and has $C^2=0$. Curves $E_{10}$, $E_1$, $E_7$ and $E_5$ of
the subgraph are orthogonal to $C$ and define Dynkin diagram
of type $D_4$. The nef cone $NEF(X)$ should be finite in
a neighbourhood of $\br_+C$. It follows that there exists
another exceptional curve $E^\prime$ of the first kind such that
$E^\prime$ together with the curves $E_{10}$, $E_1$,
$E_7$ and $E_5$ defines an extended Dynkin diagram of the type
$\widetilde{B_4}$. It follows that $X$ has exceptional curves
with the graph $H\widetilde{D_8}$. Analyzing this graph,
one can see that $X\in FPMC$, and the graph
$H\widetilde{D_8}$ gives all exceptional curves of $X$.
Here one should use some elementary facts about polyhedra with acute
angles in hyperbolic spaces (e. g. see \cite{V2}).

Similarly one can prove that all exceptional curves of $X$ have
graph $H\widetilde{A_8}$ below if curves with square $-2$
of $X$ from (***) have the extended Dynkin diagram $\widetilde{A_8}$.

This gives:
\noindent
{\it Classification of surfaces (***):
The set of all exceptional curves of the surfaces
$X$ has one of dual
graphs $H\widetilde{E_8}$, $H\widetilde{D_8}$ or $H\widetilde{A_8}$
given below. Here a black vertex $E_i$ has
$E_i^2=-2$, and a white vertex $E_i$ has $E_i^2=-1$. An edge
$E_i$, $E_j$ of the graph means that $E_i\cdot E_j=1$. General
surfaces $X$ of these type are not elliptic} (we shall show this below).

\centerline{\hbox {
\vbox{\hbox{$E_{10}$} \hbox{$\circ$}}
\hskip-10pt
\hbox to1cm{\hrulefill}
\hskip-1pt
\vbox{\hbox{$E_9$} \hbox{$\bullet$}}
\hskip-6pt
\hbox to1cm{\hrulefill}
\hskip-1pt
\vbox{\hbox{$E_8$} \hbox{$\bullet$}}
\hskip-6pt
\hbox to1cm{\hrulefill}
\hskip-1pt
\vbox{\hbox{$E_7$} \hbox{$\bullet$}}
\hskip-6pt
\hbox to1cm{\hrulefill}
\hskip-1pt
\vbox{\hbox{$E_6$} \hbox{$\bullet$}}
\hskip-6pt
\hbox to1cm{\hrulefill}
\hskip-1pt
\vbox{\hbox{$E_5$} \hbox{$\bullet$}}
\hskip-6pt
\hbox to1cm{\hrulefill}
\hskip-1pt
\hbox{\vbox to 0.8cm{\vskip0.15cm\hbox{$E_4$} \hbox{$\bullet$}
\vskip1pt
\hbox{\hskip2pt \vrule height 0.8cm width 0.5pt} \vskip-5pt\hbox{$\bullet$}
\hbox{$E_1$}}}
\hskip-6pt
\hbox to1cm{\hrulefill}
\hskip-1pt
\vbox{\hbox{$E_3$} \hbox{$\bullet$}}
\hskip-6pt
\hbox to1cm{\hrulefill}
\hskip-1pt
\vbox{\hbox{$E_2$} \hbox{$\bullet$}}
}}

\vskip2cm

\centerline{Graph $H\widetilde{E_8}$}

\vskip1cm

\centerline{\hbox {
\vbox{\hbox{$E_{11}$} \hbox{$\circ$}}
\hskip-10pt
\hbox to1cm{\hrulefill}
\hskip-1pt
\vbox{\hbox{$E_{10}$} \hbox{$\bullet$}}
\hskip-6pt
\hbox to1cm{\hrulefill}
\hskip-1pt
\hbox{\vbox to 0.8cm{\vskip0.15cm\hbox{$E_7$} \hbox{$\bullet$}
\vskip1pt
\hbox{\hskip2pt \vrule height 0.8cm width 0.5pt} \vskip-5pt\hbox{$\bullet$}
\hbox{$E_1$}}}
\hskip-6pt
\hbox to1cm{\hrulefill}
\hskip-1pt
\vbox{\hbox{$E_5$} \hbox{$\bullet$}}
\hskip-6pt
\hbox to1cm{\hrulefill}
\hskip-1pt
\vbox{\hbox{$E_3$} \hbox{$\bullet$}}
\hskip-6pt
\hbox to1cm{\hrulefill}
\hskip-1pt
\vbox{\hbox{$E_2$} \hbox{$\bullet$}}
\hskip-6pt
\hbox to1cm{\hrulefill}
\hskip-1pt
\hbox{\vbox to 0.8cm{\vskip0.15cm\hbox{$E_6$} \hbox{$\bullet$}
\vskip1pt
\hbox{\hskip2pt \vrule height 0.8cm width 0.5pt} \vskip-5pt\hbox{$\bullet$}
\hbox{$E_4$}}}
\hskip-6pt
\hbox to1cm{\hrulefill}
\hskip-1pt
\vbox{\hbox{$E_8$} \hbox{$\bullet$}}
\hskip-6pt
\hbox to1cm{\hrulefill}
\hskip-1pt
\vbox{\hbox{$E_9$} \hbox{$\circ$}}
}}

\vskip1.5cm

\centerline{Graph $H\widetilde{D_8}$}

\vskip1cm

\centerline{
\hbox to8cm {\hfill
\vbox{\hbox{$\circ$} \hbox{$E_{12}$}}
\hskip-10pt
\vbox to 1.9cm{\hbox{\hskip0.6cm$E_1$}\vfill
\hbox to1cm{\hrulefill}\vskip15pt}
\hskip-5pt
\hbox{
\vbox{\hbox{$E_4$}
\hbox{$\bullet$}
\hbox{\hskip2pt \vrule height 0.8cm width 0.5pt}
\hbox{$\bullet$}
\hbox{\hskip2pt \vrule height 0.8cm width 0.5pt}
\hbox{$\bullet$}
\hbox{$E_9$}}
\hskip-6pt
\vbox{\hbox to1cm{\hrulefill}\vskip15pt}
\hskip-1pt
\vbox to 3cm{
\hbox{\hskip-1.1cm \hbox to1.6cm{\hrulefill}}
\vfill
\hbox{$\bullet$} \hbox{$E_6$}}
\hskip-9pt
\vbox to 4.5cm{
\hbox{\hskip11pt $E_{10}$}
\hbox{\hskip11pt $\circ$}
\hbox{\hskip11pt \hskip2pt \vrule height 0.6cm width 0.5pt}
\vskip-4pt
\hbox{\hskip11pt $\bullet$}
\hbox{\hskip11pt $E_7$}
\vfill
\hbox to1.1cm{\hrulefill}\vskip15pt}
\hskip-1pt
\vbox to 3cm{
\hbox{\hskip-0.4cm \hbox to1.8cm{\hrulefill}}
\vfill
\hbox{$\bullet$} \hbox{$E_3$}}
\hskip-33pt
\vbox{\hbox to1cm{\hrulefill}\vskip15pt}
\hskip-3pt
\hbox{
\vbox{
\hbox{$E_2$}
\hbox{$\bullet$}
\hbox{\hskip2pt \vrule height 0.8cm width 0.5pt}
\hbox{$\bullet$}
\hbox{\hskip2pt \vrule height 0.8cm width 0.5pt}
\hbox{$\bullet$}
\hbox{$E_8$}
}}
\hskip-11pt
\hbox{
\vbox to 1.9cm{
\hbox{$E_5$}
\vfill
\hbox to1cm{\hrulefill}\vskip15pt
}
\hskip-2pt
\vbox{
\hbox{$\circ$} \hbox{$E_{11}$}
}
}}\hfill}
}

\vskip1cm

\centerline{Graph $H\widetilde{A_8}$}

\vskip1cm

Considering a sequence of contractions of
exceptional curves of the first kind of the surfaces $X$
(their preimages are numerated in descending
order by vertices $E_i$ of the graphs
$H\widetilde{E_8}$, $H\widetilde{D_8}$ and $H\widetilde{A_8}$),
one obtains existence and description of parameters of the
surfaces $X$. Surfaces $X$ with graphs
$H\widetilde{E_8}$, $H\widetilde{D_8}$ or $H\widetilde{A_8}$
are obtained from a line $E_1$, two different
lines $E_1$, $E_2$ and three non-collinear
lines $E_1$, $E_2$, $E_3$ respectively on a plane
${\Bbb P}^2$ by blow up of appropriate sequences of $9$ points.
They correspond to vertices of the graphs
$H\widetilde{E_8}$, $H\widetilde{D_8}$ and $H\widetilde{A_8}$
respectively in increasing order.

Let us show that general surfaces from (***) are not elliptic:
for any natural $m$ the linear system $|-mK|$
contains only the divisor $-mD$ and is
zero-dimensional. First, let us consider $m=1$. Suppose that
for a surface $X$ from (***) the linear system $|-K_X|$ is
1-dimensional. To be concrete, assume that
$X$ has the graph $H\widetilde{A_8}$ of exceptional curves.
Let $X_1$ be the surface obtained by contraction of the curve $E_{10}$.
Then $|-K_{X_1}|$ contains a nef divisor
$D_1=E_1+E_4+E_7+E_2+E_5+E_8+E_3+E_6+E_9$ with
$D_1^2=1$. By Kawamata-Viehweg vanishing
\cite{Kaw}, \cite{Vie} and Riemann-Roch theorem,
one than gets $\dim|-K_{X_1}|=1$.
Since $|-K_{X_1}|$ contains the image of
$|-K_X|$ with $\dim|-K_X|=1$, we obtain that the linear
system $|-K_{X_1}|=|-K_X|$ is one-dimensional with the base point
on the curve $E_7$. This base point is equal to the image of $E_{10}$.
Let $X^\prime$ be a surface which is blow up of $X_1$ in a
point of the curve $E_7$ which is different from
the base point of $|-K_{X_1}|$, and from points
$E_7\cap E_4$ and $E_7\cap E_2$. The surface $X^\prime$
belongs to (***), has the same graph
$H\widetilde{A_8}$ of exceptional curves, and
$\dim |-K_{X^\prime}|=0$. This shows that surfaces
(***) have one additional
parameter to elliptic surfaces $X$ from (***) with not
zero-dimensional $|-K_X|$.
Now suppose  that $|-mK|$ is not zero-dimensional, but
$|-(m-1)K|$ is zero-dimensional for some $m>0$.
Then $|-mK|$ is elliptic
pencil with the fiber $mD$ of multiplicity $m$ where $D$
has the type $\widetilde{A_8}$.
Below we shall see that for any fixed
$m\ge 1$ these elliptic surfaces have the same number of parameters
as for $m=1$. Thus, our consideration above
for $m=1$ shows that general surfaces $X$ from (***)
are not elliptic: they have one additional parameter to
elliptic surfaces from (***).

Now let us consider classification of elliptic surfaces $X$
from (*). We denote this class of surfaces by (**). It consists of:

\noindent
(**){\it Rational surfaces $X$ with $K^2=0$ having an
elliptic pencil $|-mK|$ for some $m>0$ (equivalently, $X$ is
rational minimal elliptic surface)
and such that the sum $\sum_i{r_i}=8$ for ranks $r_i$ of
reducible fibers of the pencil
(these rational minimal elliptic surfaces are called maximal).
The invariant $m$ is called index (it was first observed in \cite{Ha}).}
See \cite{D} and \cite{CD, Sect. 5.6} about these surfaces.

Classification of surfaces from (**)
can be obtained using general Ogg--Shafarevich
theory of elliptic surfaces \cite{O}, \cite{Sh3} , applied to the
special case of rational elliptic surfaces. This was
done in \cite{D} and \cite{CD}. In fact, below, we just review
these results.

Suppose that the index $m=1$.
Then the elliptic pencil $|-K_X|$ is Jacobian, it has
a section which is defined by an exceptional curve $E$ of the first
kind. For this case, there exists a sequence of contractions of
$9$ exceptional curves of the first kind $\pi:X\to {\Bbb P}^2$ such
that image of the pencil $|-K_X|$ is a pencil of plane cubics.
Condition $\sum_i{r_i}=8$ is then equivalent to
finiteness of the Mordell--Weil group defined by sections (i. e.
exceptional curves of the first kind of $X$). Using this description,
one can classify these surfaces $X$ according to Dynkin diagrams
$\Gamma_i$ of reducible fibers. It is known \cite{CD} that there
are the following and the only following possibilities for
$\Gamma_i$ with $\sum_i{r_i}=8$, and Mordell--Weyl groups,
see \cite{CD, Corollary 5.6.7}:

{\settabs 2 \columns
\+\hskip1cm {\bf Types of fibers} & {\bf Mordell--Weil groups}\cr
\+\hskip1cm $\widetilde{E}_8$ & $(1)$ \cr
\+\hskip1cm $\widetilde{D}_8$ & $\bz/2\bz$ \cr
\+\hskip1cm $\widetilde{A}_8$ & $\bz/3\bz$ \cr
\+\hskip1cm $\widetilde{E}_7+\widetilde{A}_1(\widetilde{A}_1^\ast)$ &
$\bz/2\bz$ \cr
\+\hskip1cm $\widetilde{A}_7+\widetilde{A}_1(\widetilde{A}_1^\ast)$ &
$(\bz/2\bz)^2$\cr
\+\hskip1cm $\widetilde{E}_6+\widetilde{A}_2(\widetilde{A}_2^\ast)$ &
$\bz/3\bz$ \cr
\+\hskip1cm $\widetilde{D}_5+\widetilde{A}_3$                       &
$\bz/4\bz$ \cr
\+\hskip1cm $2\widetilde{D}_4$                                      &
$(\bz/2\bz)^2$\cr
\+\hskip1cm $2\widetilde{A}_4$,                                     &
$\bz/5\bz$ \cr
\+\hskip1cm $\widetilde{D}_6+2\widetilde{A}_1(\widetilde{A}_1^\ast)$ &
$(\bz/2\bz)^2$\cr
\+\hskip1cm $\widetilde{A}_5+\widetilde{A}_1(\widetilde{A}_1^\ast)+
\widetilde{A}_2(\widetilde{A}_2^\ast)$                    &
$(\bz/3\bz)\oplus(\bz/2\bz)$\cr
\+\hskip1cm $2\widetilde{A}_3+2\widetilde{A}_1(\widetilde{A}_1^\ast)$ &
$(\bz/4\bz)\oplus (\bz/2\bz)$\cr
\+\hskip1cm $4\widetilde{A}_2(\widetilde{A}_2^\ast)$ & $(\bz/3\bz)^2$. \cr}

\smallpagebreak

It is not difficult to extend these possible diagrams of exceptional
curves with square $-2$ adding exceptional curves of the first kind.
Their number is equal to the order of the Mordell--Weil group.
E. g. for diagrams
$\widetilde{E}_8$, $\widetilde{D}_8$ and $\widetilde{A}_8$ one gets
graphs $H\widetilde{E}_8$, $H\widetilde{D}_8$ and
$H\widetilde{A}_8$ respectively given above.

Now assume that the index $m>1$.
Then the elliptic pencil $|-mK_X|$ has
a unique multiple fiber $mD$ where $D\in |-K_X|$. The Jacobian
fibration $J$ of $X$ is also a rational minimal elliptic
surface (see \cite{CD, Proposition 5.6.1}) with the same base and
the same fibers as $X$.
We have considered rational Jacobian fibrations $J$ above.
The multiple fiber $mD$ of $|-mK_X|$ fixes a point
$x_D\in D$ of order $m$ of the corresponding fiber $D$ of $J$.
By Ogg--Shafarevich theory \cite{O} and \cite{Sh3},
the triplet $(x_D\in D\subset J)$ defines
the elliptic pencil $|-mK_X|$ uniquely, any triplet
$(x_D\in D\subset J)$ (where $J$ is a Jacobian rational minimal elliptic
surface, $D$ its fiber and $x_D\in D$ its point of order $m$)
is possible and defines a rational elliptic
surface $X$ with elliptic pencil $|-mK_X|$ of index $m$,
multiple fiber $mD$ and the Jacobian fibration $J$.
See \cite{D} and \cite{CD, Ch. 5} for details.
In particular, considering of the $X$ with a multiple fiber $mD$ where
$D$ has the type $\widetilde{A_8}$ (this case has been considered above),
is equivalent to considering of triplets $(x_D\in D\subset J)$
with the reducible fiber $D$ of the type $\widetilde{A_8}$ and
the point $x_D\in D$ of order $m$. For the fixed $J$,
the set of possible $x_D\in D$ is obviously finite,
and the number of parameters of the $X$ is equal to
the number of parameters of the $J$. We have used this above.

We emphasize that the index $m$ of a surface
$X$ from (**) can be arbitrary $m\in \bn$.
It follows that the number of connected components of the moduli space
of surfaces $X$ in (**) is infinite. It follows that
{\it number of connected components of moduli space of surfaces
$X\in FPMC_{\delta_E=2,p_E=0}$ of the case 2b is infinite: it has
three connected components with non-elliptic general $X$ and
infinite number of connected components with elliptic $X$ corresponding
to infinite number of possible indexes $m\in \bn$ of $X$.}
\endremark

\smallpagebreak

It seems, nobody tried to classify surfaces
$X\in FPMC_{\delta_E=2,\,p_E}$ for $p_E\ge 1$.

\smallpagebreak

For a surface $X\in FPMC_{\rho\ge 3}$, the very important invariant
is the {\it dual graph \linebreak $\Gamma (\Exc(X))$
of the set $\Exc(X)$ of exceptional curves of $X$.} Here we mark
vertices $E\in \Exc(X)$ of this graph by the pair $(E^2,\,p_a(E))$, and
edges $(E_i, E_j)$ of the graph by $E_i\cdot E_j$ if $E_i\cdot E_j>0$.

Using considerations in the proof of Theorem 1.1, we can prove:

\proclaim{Theorem 1.2} For fixed invariants $\rho\ge 3$, $\delta_E$
and $p_E$, the set of possible graphs $\Gamma (\Exc(X))$ of
$X\in FPMC_{\rho,\delta_E,p_E}$ is finite if $K_X^2>0$.
\endproclaim

\demo{Proof} We argue as in the proof of Theorem 1.1. For the fixed
matrix $\Gamma$ (one from a finite set) we have that the hyperbolic
lattice $\NS(X)$ is an intermediate lattice
$[E_1,\dots,E_\rho]\subset \NS(X)\subset [E_1,\dots ,E_\rho]^\ast$
where the lattice $[E_1,\dots,E_\rho]$ is fixed by $\Gamma$. Thus
there exists only a finite number of possibilities for the overlattice
$[E_1,\dots,E_\rho]\subset \NS(X)$. We fix one of them. We also fix one
of finite possibilities for the canonical class
$K=K_X\in \NS(X)$. If $K^2>0$, there exists only a finite set of elements
$e\in \NS(X)$ such that $-\delta_E\le e^2<0$ and
$0\le {e^2+e\cdot K\over 2}+1\le p_E$ because the lattice $\NS(X)$ is
non-degenerate and hyperbolic.
It shows that number of possible graphs $\Gamma(\Exc(X))$ is finite.
This finishes the proof.
\enddemo

The condition $K_X^2>0$ of Theorem 1.2 is necessary.
Surfaces $X$ of the case 2b of the example 1.4.1 have
$K_X^2=0$, have infinite number of connected components of the
moduli space and infinite number of possible graphs $\Gamma(\Exc(X))$.

\head
2. Algebraic surfaces with some locally polyhedral Mori cone
\endhead

Here we want to outline some generalization of results of Sect. 1 for
more general class of surfaces (we hope to give details in
forthcoming publications). They are Algebraic Geometry
analog of reflection groups of hyperbolic
lattices of elliptic, parabolic or hyperbolic type
(see \cite{N8}, \cite{N9}, \cite{N11}, \cite{N13}---\cite{N16}).

For surfaces $X\in FPMC_{\rho\ge 3}$ the nef cone defines an
elliptic (i. e. finite and of finite volume) polyhedron
$\NEF(X)/\br_+\subset \La(X)$. The key Lemma 1.1 which we used in
the proof of Theorem 1.1, can be generalized (with some bigger
absolute constant instead of $62$) for locally
finite polyhedra
of restricted parabolic or restricted hyperbolic type
in hyperbolic spaces.
See \cite{N11}, \cite{N13}, \cite{N15} and \cite{N16}.
Thus, we shall have Theorem 1.1 for surfaces with
$\NEF(X)/\br_+$ of these types.
Below we introduce surfaces for which this is true.

Like in Sect. 1, we consider only non-singular projective algebraic
surfaces $X$ over algebraically closed field.

\definition{Definition 2.1} Let $\rho=\rho(X)\ge 3$.
We say that $X$ has {\it almost finite
polyhedral Mori cone} $\overline{NE}(X)$ if (1), (2) and (3) below
hold:

(1) There exist finite maximums:
$$
\delta_E(X)=\max_{C\in \Exc(X)}{-C^2}\ \ \text{and}\ \
p_E(X)=\max_{C\in \Exc(X)}{p_a(C)}.
$$

(2) There exists a non-zero
$r \in \NS(X)$ such that any extremal ray of $\overline{\NE}(X)$
is either generated by an exceptional curve $E\in \Exc(X)$ or
by $c\in \NS(X)\otimes \br$ such that $c^2=0$ and $c\cdot r=0$.

(3)  The set $\Exc(X)\cdot r$ is bounded:
$-R\le \Exc(X)\cdot r\le R$ for some finite $R>0$.

We remark that if $K\not\equiv 0$, one can always put $r =K$.
Then (3) follows from (1).
\enddefinition

There are plenty of surfaces $X$ with almost finite polyhedral Mori cone.
Let $Y$ be a normal projective algebraic surface with nef anticanonical
class $-K_Y$ and with at least one non Du Val singularity if
$-K_Y\equiv 0$. Then the minimal resolution $X$ of singularities of $Y$
has almost finite polyhedral Mori cone (e. g. see \cite{N9}).

One can show that surfaces $X$ with almost finite polyhedral Mori cone
have the polyhedron $\NEF(X)/\br_+\subset \La(X)$ either of elliptic
or restricted parabolic or restricted hyperbolic type. For polyhedra $\M$
of this type one can prove Lemma 1.1 with the constant $62$ replaced by
some other absolute constant. See
\cite{N11}, \cite{N13}, \cite{N15}, and \cite{N16}.
From this Lemma we get (like in proof of Theorem 1.1)

\proclaim{Theorem 2.1} For $\rho\ge 3$,
there are constants $N(\rho,\,\delta_E)$ and $N'(\rho,\,\delta_E,\,p_E)$
depending only on $(\rho,\,\delta_E)$ and $(\rho,\,\delta_E,\,p_E)$
respectively such that for any $X$ with almost finite polyhedral
Mori cone and $\rho(X)=\rho$, $\delta_E(X)=\delta_E$ and
$p_E(X)=p_E$, there exists an ample effective divisor $h$ such that
$h^2\le N(\rho,\,\delta_E)$, and if the ground field is $\bc$,
there exists a very ample divisor $h'$ such that
${h'}^2\le N'(\rho,\,\delta_E,\,p_E)$.
\endproclaim

\end